# TRAFFIC NETWORK PARTITIONING FOR HIERARCHICAL MACROSCOPIC FUNDAMENTAL DIAGRAM APPLICATIONS BASED ON FUSION OF GPS PROBE AND LOOP DETECTOR DATA

KANG AN [a], Xianbiao Hu [b*] and Xiaohong Chen [c]

[a] *Key Laboratory of Road and Traffic Engineering of the State Ministry of Education, School of Transportation Engineering, Tongji University, Shanghai 201804, China*
Email: 1310061@tongji.edu.cn

[b*] *Department of Civil, Architectural and Environmental Engineering, Missouri S&T. Rolla, MO 65401, USA.* Email: xbhu@mst.edu

[c] *Key Laboratory of Road and Traffic Engineering of the State Ministry of Education, School of Transportation Engineering, Tongji University, Shanghai 201804, China,*
Email: chenxh@tongji.edu.cn

## 1. INTRODUCTION

Macroscopic Fundamental Diagram (MFD), as a graphical method to characterize traffic state in a homogeneous road network, enables researchers and practitioners to model, monitor and control traffic at an aggregate way without the need to go into detailed link information. Various network partitioning methods were proposed in prior studies to partition a heterogeneous network into multiple homogeneous sub-networks. These approaches are mostly based on a normalized cut mechanism, which takes the traffic statistics of each link (e.g. link volume, speed or density) as input to calculate the degree of similarity between two links, and perform graph cut to divide a network into two sub-networks at each iteration when the traffic dynamics between links are dramatically different. These methods assume complete link-level traffic information over the entire network, e.g. the accurate measurement of the traffic conditions exist for every single link in the network, which makes them inapplicable when it comes to real-world setting. In a previous research by the authors [1, 2], an approach based on region growing technique and lambda-connectedness definition was proposed to partition a network while allowing certain degree of data missing ratio. The sensitivity analysis result shows traffic engineers still need to install loop detectors on at least 60% of the links in order to get good quality network partitioning and calibration result, which is still expensive to set up and maintain.

Probe vehicles with GPS devices have received much attention for its potential as a new sensing technology in the past decade [3-6]. While probe sensors can record instantaneous vehicle position and speed information, which allows us to track vehicle trajectory and reflects traffic dynamics along the traffic network over both time and space dimensions, they're mostly used to derive link-level traffic dynamics such as average travel time, delay and speed [7, 8], Having sufficient GPS probe vehicle data (PVD) coverage over the entire network to ensure a good amount of sample size on any link in the network at any time becomes a prerequisite of performing analytics and research with this approach. As a result, partitioning traffic network with GPS probe vehicle data as the sole data source becomes challenging, if practicable at all. The other main issue associated with PVD utilization in MFD related research is that previous studies usually assume a homogeneous distribution of the probe vehicles [9], i.e. an equal market penetration rate over the entire network, which is rather unrealistic.

In this paper, we propose a method which, based on fusing PVD and loop detector data (LDD), extracts the locally homogeneous subnetworks with a grid-level network approach instead of dealing with detailed link-level network. By fusing the two data sources, we take advantage of both better coverage from PVD and the full-size detection from LDD. The concept we're trying to explore and prove in this research is that, as opposed to the traditional approach of aggregating data to generate link-level traffic statistics and then perform network partition, can we approach this problem with a grid-level solution? Here, a grid is defined as the basic unit in the traffic network, which comprises multiple roadway nearby segments and can be merged to a large subnetwork or further divided into smaller grid if needed. The basic idea is that by fusing available PVD and LDD within this grid, an estimate of traffic dynamics will be computed, and then network partitioning algorithm will be further developed to perform network growing and partitioning based on the similarities between adjacent grids. With such approach, the key requirement on data size of having sufficient probe vehicles at any given link will be relaxed to having a reasonable number of probe vehicles in the vicinity to capture the local traffic dynamics, which becomes more feasible and practicable.

Another natural merit of this modeling approach is that the penetration rate of the probe vehicles, which is a key variable in measuring and calibrating MFD [10] and is calculated by the ratio of number of probe vehicles to the total traffic counts from loop detectors, can be computed for each grid at different time slots, instead of assuming a single value for the entire network. It's also worth noting that for the data from loop detectors, only the volume is used but not the speed or density information, which makes this proposed model generically applicable to not only freeway but also arterial network, regardless of which part of road segment are they installed.

## 2. METHODOLOGY

### 2.1. Network Density Estimation

We start from a simple network where there's only one loop detector available, represented by Figure 1(a). A total of N vehicles travel in this network and all of them follow a homogenous path that passes loop detector, with a path represented in solid line, whereas all untraveled links are drawn in dash lines. Among them, trajectory data of $P$ vehicles are available ($P \leq N$).

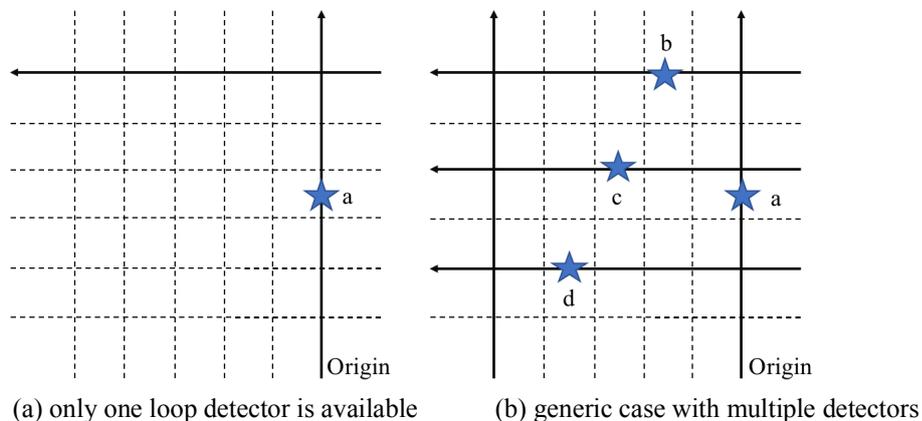

(a) only one loop detector is available     (b) generic case with multiple detectors
Figure 1 network illustration

Let's use $d_i$ and $t_i$ ($i = 1,2,...S$) to denote the distance traveled and travel time of vehicle $i$, respectively, within the space-time window defined by link $l$ and time interval $[t, t + \Delta t]$.

Assuming the probe vehicle is randomly sampled from the entire vehicle set and is able to represent the population, the total travel time by all vehicles can be calculated by

$$T = \sum_{i=1}^{P} t_i * N/P = \sum_{i=1}^{P} t_i /p \qquad \text{Eq.(1)}$$

Where $p$ is the sample rate, or sometimes called market penetration rate, and is calculated by

$$p = P/N \qquad \text{Eq.(2)}$$

Similarly, we will have the following for total distance traveled by all vehicles

$$D = \sum_{i=1}^{P} d_i * N/P = \sum_{i=1}^{P} d_i /p \qquad \text{Eq.(3)}$$

With this, density of this network can be computed with Edie's generalized definitions, i.e.

$$K = T/l\Delta t = \sum_{i=1}^{P} t_i /pl\Delta t \qquad \text{Eq.(4)}$$

Similarly, average volume and speed of this network can be computed

$$Q = D/l\Delta t = \sum_{i=1}^{P} d_i /p\Delta t \qquad \text{Eq.(5)}$$

$$V = Q/K = \frac{\sum_{i=1}^{P} d_i}{p\Delta t} / \sum_{i=1}^{P} t_i /pl\Delta t = \sum_{i=1}^{P} d_i / \sum_{i=1}^{P} t_i \qquad \text{Eq.(6)}$$

Next, we move to a generic network where multiple loop detectors are available, and vehicles are traveling in the network following multiple heterogeneous paths, represented by Figure 1(b). In this particular example, detectors $a$ and $b$ are located on the same path, whereas detectors $a$ and $c$ are located on another same path, and detector $c$ is on another path by itself. With the existence of multiple detectors, the calculation of market penetration rates becomes not only spatially dependent, but more importantly, the interactions of multiple sensors needs to be considered. In other words, as one vehicle may travel through multiple detectors, it may contribute to the calculation of market penetration rate and network density for multiple times, and such redundancy brings bias to the final result. For example, if a vehicle travel through both sensors $a$ and $b$, the trajectory data of this vehicle which essentially represent the behavior of driver behind the steering wheel, will be used twice, as opposed to only one time for another vehicle with same origin and destination locations but follows a different path via detector $d$.

To address this issue, we introduced the following variables

$M$: total number of loop detectors
$Veh_i$: $i^{\text{th}}$ vehicle
$\emptyset_{ij}$: vehicle-detector relationship.
$w_{ij}$: weight of detector $j$ for vehicle $i$

We will have:

$$\emptyset_{ij} = \begin{cases} 1, \text{if vehicle } i \text{ goes through loop detector } j \\ 0, \text{otherwise} \end{cases} \qquad \text{Eq.(7)}$$

and

$$w_{ij} = \emptyset_{ij} / \sum_{k=1}^{M} \emptyset_{ik} \qquad \text{Eq.(8)}$$

Eq. (2) will turn into

$$p_j = \sum_{i=1}^{P} \emptyset_{ij} /N_j \qquad \text{Eq.(9)}$$

$$p = \sum_{j=1}^{M} \sum_{i=1}^{P} \emptyset_{ij} / \sum_{j=1}^{M} N_j \qquad \text{Eq.(10)}$$

Eq. (4) will turn into

$$K = \sum_{j=1}^{M} \sum_{i=1}^{P} w_{ij} * t_i /p_j l\Delta t$$
$$= \sum_{j=1}^{M} \sum_{i=1}^{P} \emptyset_{ij} * t_i * N_j / \sum_{k=1}^{M} \emptyset_{ik} \sum_{i=1}^{P} \emptyset_{ik} \, l\Delta t \qquad \text{Eq.(11)}$$

### 2.2. Pseudocode

---

**I**: Initial grid definition

For each $TAZ$:

    N= (int) size of $TAZ$/minimum size

    Partition each TAZ to N grids, each with minimum size

**II**: grid attributes calculation

For each partition

    Calculate market penetration rate $p$ via Eq. (10)

    Calculate traffic density $K$ via Eq. (11)

**III**. network conversion

Build new nodes, links and connectivity

$C(i,j) = 1/e^{(k_i-k_j)^2}$  #potential function

**IV**. partition growing

    $\bar{S}=\emptyset$; S={N}

    $r = 0$; #region ID

    While $S!=\{\}$:

        $r = r+1$

        $seed_r = \underset{m \in S}{argmax}\ C_m$

        $\bar{S}(r) = \bar{S}(r) \cup \{seed_r\}, S = S \setminus \{seed_r\}$

        *While $\bar{S}(r)!=\{\}$:*

            $i = \underset{m \in \bar{S}}{argmax}\ C_m$

            For each $j \in A(i)$:

                If $C(i,j) \geq \lambda$:

                    $\bar{S}(r) = \bar{S}(r) \cup \{j\}, S = S \setminus \{j\}$

                    $C_j = \min(C_j, C(i,j))$

                Else:

                    Pass

    Output: $\bar{S}$ as the partitioning result

**V**: integrative region growing

Repeat steps II~IV until no more partitions can be further merged

Smooth network partition boundaries

---

## 3. PRELIMINARY ANALYSIS

Figure 2 shows the network partition result. While further analysis is under way, preliminary analysis shows that compared with previous approach that requires a single data source of 60% loop detector data, the proposed new approach is able to achieve comparable result with a scenario of data availability as low as 15%.

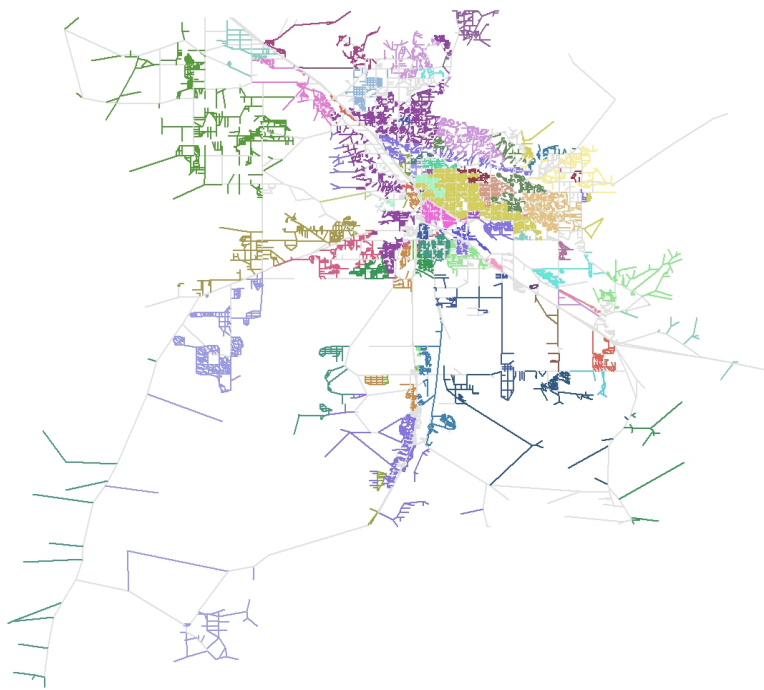

Figure 2 network partition result

Other numeric analysis in progress include analysis of minimum grid size, heuristics of grid definition, sensitive analysis with regard to data availabilities and other aspects.